\documentclass[12pt]{article}

\usepackage{amsmath,amssymb,amsfonts,amsthm}
\def\Q{\mathbb{Q}}

\def\F5{\Bbb F_5}

\def\C{\Bbb C}

\def\e2{{\eta_2}}
\def\ei{{\eta_i}}

\def\no{\noindent}

\def\p{\pi}

\def\ti{\times}

\def\ot{\otimes}
\def\part{\partial}

\begin{document} 

\title{Murre's conjectures for certain product varieties}
\author{Kenichiro Kimura}
\date{}
\maketitle

\begin{abstract}
  We consider Murre's conjectures on Chow groups for a fourfold which is a
product of two curves and a surface. We give a result which concerns Conjecture
D:the kernel of a certain projector is equal to the homologically trivial part
of the Chow group. We also give a proof of Conjecture B for a product of two
surfaces.  MSC number: 14C25

\end{abstract}

\section{Introduction}

Let $X$ be a smooth projective variety over $\C$ of dimension 
$d$. Let $\Delta
\subset X\ti X$ be the diagonal. There is a cohomology class 
$cl(\Delta)\in H^{2d}(X\ti X).$

\no In this paper we use Betti
cohomology with rational coefficients. There is the K\"unneth decomposition

\hfil $H^{2d}(X\ti X)\simeq \bigoplus_{i=0}^{2d} H^{2d-i}(X)\ot H^i(X).$\hfil

\no We write $cl(\Delta)=\sum_{i=0}^{2d}\pi_i^{hom}$ according to 
this decomposition. Here $\pi_i^{hom}\in H^{2d-i}(X)\ot H^i(X).$
If the K\"unneth conjecture is true, then each $\pi_i^{hom}$ 
is an algebraic cycle.

\no Murre(\cite{Mu},\cite{Mu2})
 formulated the following conjecture. For an abelian group $M$, 
we write $M_\Q=M\ot\Q.$ 

\begin{enumerate}
\renewcommand{\labelenumi}{(\Alph{enumi})}
\item  The $\pi_i^{hom}$ lift to a set of orthogonal projectors $\pi_i$
in $CH^d(X\ti X)_\Q$ which satisfy the equality

 \[ \sum_{i=0}^{2d}\pi_i =\Delta.\]
 
\item The correspondences $\pi_0,\cdots ,\pi_{j-1},\,
\pi_{2j+1},\cdots , \pi_{2d}$ act as zero on $CH^j(X)_\Q.$

\item Let $F^\nu CH^j(X)=Ker \pi_{2j}\cap Ker \pi_{2j-1}
\cdots \cap Ker \pi_{2j-\nu+1}.$  Then the filtration $F^\cdot$ is independent
of the choice of $\pi_i$.

\item $F^1CH^j(X)_\Q=CH^j(X)_{hom,\,\Q}.$

\end{enumerate}

\no It is shown by Jannsen(\cite{Ja}) that 
this conjecture of Murre is equivalent
to Beilinson's conjectures on the filtrarion on Chow groups.

\no There are not yet many evidences for this conjecture. For a projective 
smooth curve $C$ and a closed point $p$ on $C$, 
set $\p_0=p\ti C,\,\,\p_2=C\ti p$ and $\p_1=\Delta-\p_0-\p_2.$
Then Conjectures (A), (B) and (D) are true for these projectors. 
For a projective smooth surface  Murre(\cite{Mu}) constructed
a set of  projectors $\p_0,\cdots , \p_4$ for which Conjectures (A), (B)
and (D) are true. About Conjecture (C) he proved that the filtration 
on Chow groups given by these projectors is a natural one in the following sense(Theorem 3 in \cite{Mu}):
\begin{itemize}

\item $F^1(CH^1(S)_\Q)=Ker(\p_2)=Pic^0(S)_\Q.$

\item $F^1(CH^2(S)_\Q)=CH^2(S)_{hom,\Q}.$ 
$F^2(CH^2(S)_\Q)=Ker(\p_3)=Ker(alb: CH^2(S)_{hom,\Q} \to Alb(S)_\Q).$

\end{itemize}

Conjecture (A) is also true for 
abelian varieties (Shermenev \cite{Sh}, Deninger-Murre \cite{DM}), 
hypersurfaces (easy), certain class of 
threefolds (del Angel-M\"uller-Stach \cite{deM},\cite{deM2}), and some modular varieties 
(Gordon-Murre \cite{GM}, Gordon-Hanamura-Murre \cite{GHM}, \cite{GHM2},
Miller-M\"uller-Stach-Wortmann-Yang-Zuo \cite{Pic}). 

\no Note that if Conjecture (A) is true for varieties $X$ and $Y$,
then it is also true for $X\ti Y$. One can put 
${\pi_i}_{X\ti Y}=\sum_{p+q=i}{\pi_p}_X\ti {\pi_q}_Y.$

\no In \cite{Mu2} Murre proves that Conjectures (B) and (D) are true
for a product of a curve and a surface for this product Chow-K\"unneth
decomposition.

Recently Murre(\cite{KMP}) proved the validity of Conjecture (B) 
and some part of Conjecture (D) for a product of two surfaces. 
More precisely, Murre proved that Conjecture (D) is true 
for a product $S_1\ti S_2$
of two smooth projective surfaces except the following part:

\no The projector ${\pi_2}_{S_1}\ti {\pi_2}_{S_2}$ act as zero
on $CH^2(S_1\ti S_2)_{hom,\,\Q}.$

\no If this is true for the case of a self-product $S_1=S_2$ of a surface, 
then 
Bloch's conjecture ($p_g=0 \Rightarrow$ albanese map is injective)
for $S_1$ is true.  If one assumes that the Chow group of $S_1$ is 
finite dimensional in the sense of Kimura(\cite{Ki}), then for an element
$z\in CH^2(S_1\ti S_1)_{hom,\Q}$ one has the equality

\[(\p_2\ti \p_2(z))^n=0\] where ${}^n$ means the power as a correspondence
and $n$ is determined by the second Betti number of $S_1$.

\no In this paper we consider Conjecture (D) for the case where
$X$ is a product of two curves and a surface $C_1\ti C_2\ti
S.$ In this case the most
crucial part is to show that ${\pi_1}_{C_1}\ti {\pi_1}_{C_2}\ti 
{\pi_2}_S$ act as zero on $CH^2(X)_{hom, \Q}.$ Here the projectors
${\pi_1}_{C_i}$ for $i=1$ and $2$ are defined as above
and we refer the reader to \cite{Mu} for the definition of the 
projector ${\p_2}_S$. Our original aim was to show 
that if the cohomology $H^1(C_1)\ot H^1(C_2)\ot H^2(S)$ has no
non-zero Hodge cycle, then ${\pi_1}_{C_1}\ti {\pi_1}_{C_2}\ti 
{\pi_2}_S$ kills all the codimension 2 cycles on $X$. 
We could not completely solve the problem, so instead we studied
what kind of cycles are killed by ${\pi_1}_{C_1}\ti {\pi_1}_{C_2}
\ti {\pi_2}_S$. It seems that under certain assumptions
on $X$, ``generic'' cycles are killed by this
projector (Theorem 2.1). This is the main result of this paper. 

We also give a proof of the essential part of
Conjecture (B) for a product of two 
surfaces. Our proof is similar to that of Murre in that we make
essential use of the properties of the Chow-K\"unneth projectors
for surfaces constructed by Murre. However there are still some
differences so we decided to include our proof here. 

This paper is organized as follows. In Section two we prove 
our main result about Conjecture (D). Section  three is devoted to a proof
of Conjecture (B) for a product of two surfaces.

The author expresses his deep gratitude to Jacob Murre 
for his patience in reading an earlier version of this paper and for
valuable comments and
encouragement.

\section{The main result.}

Let $C_1$ and $C_2$ be smooth projective curves over $\C$. Let $S$
be a smooth projective surface over $\C$.
We assume that these varieties are sufficiently general so that they
satisfy the following conditions:
\begin{itemize}

\item $NS(S)\ot \Q=\Q H$ where $H$ is a hyperplane section of $S$.

\item The cohomology groups 
$H^1(C_1)\ot H^1(C_2),\,H^1(C_2)\ot H^1(S)$ and
$H^1(C_1)\ot H^1(S)$ have no non-zero Hodge cycle. 

\end{itemize}

\no  Let $X=C_1\ti C_2\ti S.$

\no Let $Z$ be a closed irreducible subvariety of $X$ of codimension 2. Consider
the following conditions for $Z$.

\begin{enumerate}

\item $pr_{12}(Z)\subset C_1\ti C_2$ has dimension $\leq 1$.

\item $pr_3(Z)\subset S$ has dimension $\leq 1$.

\item $pr_{12}:Z\to C_1\ti C_2$ and $pr_3: Z\to  S$ are surjective and 
$Z$ is a Cartier divisor either of $C_1\ti pr_{23}(Z)$ or of $C_2\ti
pr_{13}(Z)$.

\end{enumerate}

\no Here we denote by $pr_*$ the various projections from $X=C_1\ti C_2\ti S.$
About the condition 3, if we assume that the projection 
$pr_{23}:Z\to pr_{23}(Z)$
is flat, then one can show that $Z$ is a Cartier divisor on 
$C_1\ti pr_{23}(Z)$(Lemma 2.2).

\newtheorem{theo}{Theorem}[section]
\begin{theo}[]
Let $C_1$ and $C_2$ be a projective smooth curves over $\C$ and let
$S$ be a projective smooth surface over $\C$. 

\no Let $X=C_1\ti C_2\ti S.$ Assume that these varieties satisfy the 
following conditions:
\begin{itemize}

\item $NS(S)\ot \Q=\Q H$ where $H$ is a hyperplane section of $S$.

\item The cohomology groups 
$H^1(C_1)\ot H^1(C_2),\,H^1(C_2)\ot H^1(S)$ and
$H^1(C_1)\ot H^1(S)$ have no non-zero Hodge cycle. 

\end{itemize}
Assume that the surface $Z$ satisfies one of the conditions 1,2 and 3 above.
Then the Chow-K\"unneth projector ${\pi_1}_{C_1}\ti {\pi_1}_{C_2}
\ti {\pi_2}_S$ kills $Z$ in $CH^2(X)_\Q$.
\end{theo}
{\it Proof.} Assume that the condition 1 holds for $Z$. 
Note that we have a factorization

\[{\pi_1}_{C_1}\ti {\pi_1}_{C_2}\ti {\pi_2}_S
=({\pi_1}_{C_1}\ti {\pi_1}_{C_2}\ti id_S)\circ (id_{C_1\ti C_2}
\ti {\pi_2}_S)\] and they commute.
We write 
$C=pr_{12}(Z)\subset C_1\ti C_2.$ Let $\eta_C\overset{j}{\hookrightarrow}
C$ be the generic point of $C$. We apply the projector
$id_{C_1\ti C_2}\ti (\pi_2)_S$ on $Z$ as a cycle on $C\ti S.$
We have the equality
\[ (j\ti id_S)^*(id_C\ti {\pi_2}_S)(Z)=(\eta_C\ti {\pi_2}_S)((j\ti id_S)^*Z).\]
We write $(j\ti id_S)^*Z=Z_\eta.$ 

\no Since $Z_\eta$ is algebraically equivalent to a cycle $\eta_C\ti E$ on 
the surface $\eta_C\ti S$  where $E$ is a divisor on $S$ 
defined over the base field
$\C$, we see that 

\hfil $(\eta_C\ti {\pi_2}_S)Z_\eta=(\eta_C\ti {\pi_2}_S)(\eta_C\ti E).$\hfil

\no Here we use that ${\p_2}_S(Pic^0(S)_\Q)=0.$ 
By taking the closure of this equality in $C\ti S$, it follows that
\[ (id_C\ti {\pi_2}_S)(Z)=C\ti {\pi_2}_S(E)+\sum_t p_t\ti S\]
where for each $t$ $p_t$ is a closed point on $C$. Applying 
$id_C\ti {\pi_2}_S$ on both sides of the equality kills $p_t\ti S$
because by Conjecture (B) for $S$ ${\pi_2}_S(S)=0.$

\no Then we apply ${\pi_1}_{C_1}\ti {\pi_1}_{C_2}\ti id_S$ on both 
sides of the equality. Since the cohomology $H^1(C_1)\ot H^1(C_2)$
has no non-zero Hodge cycle, it follows that 

\hfil ${\pi_1}_{C_1}\ti {\pi_1}_{C_2}(C)\in Pic^0(C_1\ti C_2).$\hfil

\no Since $Pic^0(C_1\ti C_2)=pr_1^*Pic^0(C_1)+pr_2^*Pic^0(C_2)$ 
by applying  ${\pi_1}_{C_1}\ti {\pi_1}_{C_2}$ again on 
${\pi_1}_{C_1}\ti {\pi_1}_{C_2}(C)$ we see that it is zero because
${\pi_1}_{C_i}CH^0(C_i)=0$ for $i=1$ and $2$.

\no The proof is similar if we assume that the condition 2 holds
for $Z$.

\no Next we assume that the condition 3 holds for $Z$. 
Assume that $Z$ is a Cartier divisor
on $C_1\ti pr_{23}(Z).$ 

\newtheorem{lem}{Lemma}[section]
\begin{lem} 
The subvariety $pr_{23}(Z)\subset C_2\ti S$ is an ample divisor. 
\end{lem}
{\it Proof}. By the assumtions on $C_2$ and $S$, we see that

\hfil $NS(C_2\ti S)\ot \Q=\Q(pt\ti S)\oplus \Q(C_2\ti H)$. \hfil

\no We denote $D_1=pt\ti S$ and $D_2=C_2\ti H$. Write $aD_1+bD_2$ for
the class
of $pr_{23}(Z)$ in $NS(C_2\ti S)\ot\Q.$ We see that 
$a=(C_2\ti pt,\,pr_{23}(Z))>0$ and $b=\frac{(pt\ti H,\,pr_{23}(Z))}{
(H,H)}>0$. Here $(*,*)$ denotes intersection number. So it follows 
that $pr_{23}(Z)-aD_1-bD_2\in Pic^0(C_2\ti S)\simeq 
Pic^0(C_2)\oplus Pic^0(S).$ So there are divisors $d_1\in Pic(C_2)$
and $d_2\in Pic(S)$ such that 
\hfil $pr_{23}(Z)=pr_2^*d_1+pr_3^*d_2$ in $Pic(C_2\ti S).$\hfil
By Nakai's criterion $d_2$ is an ample divisor on $S$ and $d_1$ is ample
on $C_2$.\hfill \qed

\no By Lemma 2.1 it follows that $H^1(pr_{23}(Z))\simeq H^1(C_2\ti S)
\simeq H^1(C_2)\oplus H^1(S).$ 
Since $Z$ is a Cartier divisor on $C_1\ti pr_{23}(Z)$, we can consider
its cohomology class 

\no $cl(Z)\in H^2(C_1\ti pr_{23}(Z))\simeq
H^2(C_1)\oplus H^1(C_1)\ot H^1(pr_{23}(Z))\oplus H^2(pr_{23}(Z))$
which is the class associated to the line bundle ${\cal O}(Z).$

\no We write $cl(Z)=c_1+c_2+c_3$ according to this decomposition.

\no Let $f:{\bf  S}\to pr_{23}(Z)$ be a resolution of singularity.
Since $Z$ is not contained in the singular locus of $C_2\ti pr_{23}(Z)$
we can take the pullback $(id\ti f)^*Z$ in $C_2\ti {\bf S}$ with the 
associated line bundle $(id\ti f)^*{\cal O}(Z).$ This pullback is
compatible with the pullback on cohomology.

\no The class $(id_{C_1}\ti f)^*c_2\in
H^1(C_1)\ot f^*(H^1(C_2)\oplus H^1(S))$ is a Hodge cycle
and by the assumption
the cohomology $H^1(C_1)\ot (H^1(C_2)\oplus H^1(S))$
has  no non-zero Hodge cycle. So we see that  $(id_{C_1}\ti f)^*c_2=0$. 

\no So there are divisors $d_1\in Pic(C_1)$ 
and $d_2\in Pic(\bf S)$ such that in 
 $Pic(C_1\ti\bf  S)$ there is an equality

\[(id_{C_1}\ti f)^*Z=d_1\ti {\bf  S}+C_1\ti d_2. \] 
 Pushing down to $C_1\ti pr_{23}(Z)$ by the map $id_{C_1}
 \ti f$ we have an equality
\[Z=d_1\ti pr_{23}(Z) +C_1\ti f_*(d_2) \] 
in $A_2(C_1\ti pr_{23}(Z))$.

\no Once $Z$ is of this form, one can see that Chow-K\"unneth projector
${\pi_1}_{C_1}\ti {\pi_1}_{C_2} \ti {\pi_2}_S$ kills $Z$ in $CH^2(X)$
because by Conjcture (B) for $C_2\ti S$ (\cite{Mu2}) 
${\pi_1}_{C_2} \ti {\pi_2}_S$ kills $ pr_{23}(Z)$
and ${\pi_1}_{C_1}$ kills $C_1$. \hfil \qed

\begin{lem}
If the projection $pr_{23}:Z\to pr_{23}(Z)$ is flat, then 
$Z$ is a Cartier divisor on $C_1\ti pr_{23}(Z)$.
\end{lem}

{\it Proof.} Let $I_Z$ be the ideal sheaf of $Z$ in $C_1\ti pr_{23}(Z).$
 For any point 
$x\in pr_{23}(Z)$, Let $\{z_i\}_i$ be the set of closed points 
on the fiber $Z \ti_{pr_{23}(Z)}{\rm Spec}\kappa(x).$ The image of $I_Z$
in the local ring 
${\cal O}_{C_1\ti_{\C}{\rm Spec}\kappa(x),\,z_i}$
is a principal ideal $(f_i).$ For each $i$ take a local section 
$\tilde{f}_i\in I_Z$ which has the image $f_i$ 
in ${\cal O} _{C_1\ti_{\C}{\rm Spec}\kappa(x),\,z_i}.$

\no For a sufficiently small neighborhood $U$ of $x$ in $pr_{23}(Z)$ we can 
consider a Cartier divisor $D$ on $C_1\ti U$ which is defined by the equation
$\tilde{f}_i$ in a neighborhood of $z_i$. Let $K$ be the kernel of natural
surjection ${\cal O}_D\to {\cal O}_Z.$ 

\no Let $\phi_D$ be the function on the set of points on $U$ defined by
\[\phi_D(y)={\rm dim}_{\kappa(y)}{\cal O}_D\ot_{{\cal O}_U}\kappa(y).\]
It is an upper semicontinuous function on $U$. So there is an neighborhood
$U'\subset U$ of $x$ such that for any $y\in U'$, one has

\hfil$\phi_D(y)\leq \phi_D(x).$ \hfil

\no On the other hand, 
${\rm dim}_{\kappa(y)}{\cal O}_Z\ot_{{\cal O}_U} \kappa(y)$
is a constant function since $Z$ is flat over $pr_{23}(Z).$ Also note that
$\phi_D(y)\geq {\rm dim}_{\kappa(y)}{\cal O}_Z\ot_{{\cal O}_U} \kappa(y)$
on $U'.$ 
Since $\phi_D(x)={\rm dim}_{\kappa(x)}{\cal O}_Z\ot_{{\cal O}_U} \kappa(x)$
it follows that 
\[\phi_D(y)= {\rm dim}_{\kappa(y)}{\cal O}_Z\ot_{{\cal O}_U} \kappa(y)\]
on $U'$. As ${\cal O}_Z$ is a flat ${\cal O}_U$ module, it follows that
\hfil $K\ot_{{\cal O}_U} \kappa(y)=0$\hfil
for any point $y\in U'.$ Hence $K=0.$ \hfill \qed

\section{A proof of Conjecture (B) for a product of two surfaces.}

In this section we give a proof of the essential part of Conjecture (B)
for a product of two surfaces.

\no Let $S_1$ and $S_2$ be projective smooth surfaces over $\C$ and
let $X=S_1\ti S_2$. For each $S_i$ there is a Chow-K\"unneth decomposition
${\p_0}_{S_i},\cdots ,
{\p_4}_{S_i}$ of the diagonal consturcted by Murre(\cite{Mu}).
 They have the following properties:

$\p_4,\,\p_3$ and $\p_0$ act as 0 on $CH^1(S_i)_\Q.$
$F^1CH^1(S_i)_\Q=Ker(\pi_2)=CH^1(S_i)_{hom, \Q}$. $F^2CH^1(S_i)_\Q
=Ker(\p_1|_{F^1})=0$.

$\p_0$ and $\p_1$ act as 0 on $CH^2(S_i)_\Q.$ 
$F^1CH^2(S_i)_\Q=Ker(\p_4)=CH^2(S_i)_{hom, \Q}.$
$F^2CH^2(S_i)_\Q=Ker(\p_3|_{F^1})
=Ker(alb:\,CH^2(S_i)_{hom,\Q}\to Alb(S_i)\ot\Q).$
$F^3CH^2(S_i)_\Q=Ker(\p_2|_{F^2})=0.$

There is a Chow-K\"unneth decomposition for $X$ given by the product 
of those for $S_i$. 

\no Murre has proven Conjecture (B) for $X$. Here we give another proof
of the essential part of his result.
 \begin{theo}[]
The Chow-K\"unneth projectors ${\p_3}_{S_1}\ti {\p_3}_{S_2}$ 
and ${\p_3}_{S_1}\ti {\p_2}_{S_2}$ act as 
zero on $CH^2(X)_\Q.$
 \end{theo}

{\it Proof.} 
Let $Z$ be an element of $CH^2(X).$
Let $\eta_i
\overset{j_i}{\hookrightarrow}
 S_i$ be the generic point of $S_i$ for $i=1,\,2$ and $Z_{\ei}$ be the generic
 fiber of $Z$.

\no The case of ${\p_3}_{S_1}\ti {\p_3}_{S_2}.$ 
 $(id_{S_1}\ti j_2)^*({\p_3}_{S_1}\ti id_{S_2}) (Z)=
 {\p_3}\ti \e2 ((id_{S_1}\ti j_2)^*Z).$ We Write ${\p_3}\ti \e2=
 {\p_3}_{\e2}$ and $(id_{S_1}\ti j_2)^*Z=Z_{\e2}.$
For $p=1$ and 2 let $C_p\overset{i_p}{\hookrightarrow}
S_p$ be a smooth hyperplane section defined over the base field 
$\C$.  Then by Lemma 2.3
of \cite{Mu}, 
${i_p}_*:\,Jac(C_p)\to Alb(S_p)$ is a surjection. So it follows that 
${i_1}_*:\,Jac(C_1)(\e2)_\Q\to Alb(S_1)(\e2)_\Q$ is also surjective. Let
$d$ be the degree of $Z_\e2$ and let $e_1$ be a closed point on $S_1$
which is rational over the base field $\C$. 
Then $Z_\e2-d(e_1)\in CH^2({S_1}_\e2)_{hom, \Q}$
and so there is a cycle $D\in Pic^0C_1(\e2)_\Q$ such that $alb(Z_\e2-d(e_1))
={i_1}_*(D).$ Let $\bar{D}$ be the closure of $D$ in $X$. Since $D$ 
is supported on $C_1\ti \eta_2$, $\bar{D}$ is supported on 
$\overline{C_1\ti \eta_2}=C_1\ti S_2$.

\no Since $Ker \p_3=Ker(alb)$, we have the equality
$$(id_{S_1}\ti j_2)^*({\p_3}_{S_1}\ti id_{S_2})(Z -d(e_1)\ti S_2-\bar{D})=
{{\p_3}_{S_1}}_\e2(Z_\e2-d(e_1)-{i_1}_*D)=0.$$ 
So it follows that 

\hfil $({\p_3}_{S_1}\ti id_{S_2})(Z)=({\p_3}_{S_1}\ti id_{S_2})
(\bar{D})+d{\p_3}_{S_1}(e_1)\ti S_2+\sum_k D_k$ \hfil

\no where for each $k$ $D_k$ is supported on 
$S_1\ti Y_k$ for an irreducible curve
$Y_k$. We apply the projector 
${\p_3}_{S_1}\ti id_{S_2}$ again on both sides of the equality.
We apply ${\p_3}_{S_1}\ti id_{S_2}$ on each $D_k$ as a cycle 
on $S_1\ti Y_k$. Let $\eta_Y\overset{j_Y}{\hookrightarrow} 
Y_k$ be the generic point of $Y_k$. We have the equality 

\hfil $(id_{S_1}\ti j_Y)^*({\p_3}_{S_1}\ti id_{Y_k})(D_k)=
({\p_3}_{S_1}\ti \eta_Y)((id_{S_1}\ti j_Y)^*D_k).$\hfil

\no Since $(id_{S_1}\ti j_Y)^*D_k$ is a divisor on the surface 
$S_1\ti \eta_Y$, from Conjecture (B) for $S_1$ it follows that 

\hfil $({\p_3}_{S_1}\ti \eta_Y)((id_{S_1}\ti j_Y)^*D_k)=0.$\hfil

\no By taking closure of this equality in $S_1\ti Y_k$ we have the equality

\hfil $({\p_3}_{S_1}\ti id_{Y_k})(D_k)=\sum_i S_1\ti p_i$ \hfil 

\no where for each $i$ $p_i$ is a closed point on $Y_k$. Applying
${\p_3}_{S_1}\ti id_{S_2}$ again on both sides of the equality
it follows that $({\p_3}_{S_1}\ti id_{S_2})(S_1\ti p_i)=0$ since
by Conjecture (B) for $S_1$ ${\p_3}_{S_1}(S_1)=0.$

\no Next we apply $id_{S_1}\ti {\p_3}_{S_2}$ on both sides of the equality. 
By Conjecture (B) for $S_2$ we see that

\hfil $(id_{S_1}\ti {\p_3}_{S_2})(d{\p_3}_{S_1}(e_1)\ti S_2)
=d{\p_3}_{S_1}(e_1)\ti {\p_3}_{S_2}(S_2)=0.$\hfil 

\no Since
$$({\p_3}_{S_1}\ti id_{S_2})(id_{S_1}\ti {\p_3}_{S_2})=(id_{S_1}\ti {\p_3}_{S_2})({\p_3}_{S_1}\ti id_{S_2})$$
it follows that 
$$(id_{S_1}\ti {\p_3}_{S_2})({\p_3}_{S_1}\ti id_{S_2})(\bar{D})
=({\p_3}_{S_1}\ti id_{S_2})(id_{S_1}\ti {\p_3}_{S_2})(\bar{D})
=({\p_3}_{S_1}\ti id_{S_2})(\sum_l p_l\ti S_2)$$ for a set of closed points $p_l$ on $S_1$. Here we apply $id_{S_1}\ti {\p_3}_{S_2}$ on $\bar{D}$ as a cycle
on $C_1\ti S_2$.
So we are reduced to the case where each component of 
$Z$ is of the form $pt\ti S_2$
 for a closed point $pt.$ We can see that the projector
 ${\p_3}_{S_1}\ti {\p_3}_{S_2}$ kills $pt\ti S_2$ bacause 
by Conjecture (B) for surfaces ${\p_3}_{S_i}(S_i)=0$ for $i=1$ and $2$.

\no{\bf Remark.} Murre pointed out that there is a simpler argument than the 
one above.
We use the equality 
\[{\p_3}_{S_1}\ti id_{S_2}(Z)=Z\circ {}^t{\p_3}_{S_1}
=Z\circ {\p_1}_{S_1}\]
where $\circ$ is composition as correspondences
and ${}^t$ is transpose. By construction of $\p_1$ there is a curve $C$ on $S_1$ such that 
${\p_1}_{S_1}$ is supported on $C\ti S_1$(cf. (ii) of Proposition 2.1 in 
\cite{KMP}). So one 
can immediately conclude that 
${\p_3}_{S_1}\ti id_{S_2}(Z)$ is supported on $C\ti S_2$.

$\quad$

\no The case of ${\p_3}_{S_1}\ti {\p_2}_{S_2}.$ We use the factorization
${\p_3}_{S_1}\ti {\p_2}_{S_2}=(id_{S_1}\ti {\p_2}_{S_2})({\p_3}_{S_1}\ti id_{S_2})$. We have the equality
$$({\p_3}_{S_1}\ti id_{S_2})(Z)=({\p_3}_{S_1}\ti id_{S_2})(\bar{D})+d{\p_3}_{S_1}(e_1)\ti S_2+\sum_k D_k$$
where for each $k$ $D_k$ is supported on $S_1\ti Y_k$ for an irreducible curve
$Y_k$ and $\bar{D}$ is supported on $C_1\ti S_2$. 
The $D_k$ part can be treated as above. 
Then we apply $id_{S_1}\ti {\p_2}_{S_2}$
on both sides of the equality. By Conjecture (B) for $S_2$ it follows that

\hfil $(id_{S_1}\ti {\p_2}_{S_2})(d{\p_3}_{S_1}(e_1)\ti S_2)=
d{\p_3}_{S_1}(e_1)\ti {\p_2}_{S_2}(S_2)=0.$\hfil

\no By using the equality
\[({\p_3}_{S_1}\ti id_{S_2})(id_{S_1}\ti {\p_2}_{S_2})=(id_{S_1}\ti {\p_2}_{S_2})({\p_3}_{S_1}\ti id_{S_2})\] we have 
$$(id_{S_1}\ti {\p_2}_{S_2})({\p_3}_{S_1}\ti id_{S_2})(\bar{D})=
({\p_3}_{S_1}\ti id_{S_2})(id_{S_1}\ti {\p_2}_{S_2})(\bar{D}).$$
Let $\eta_{C_1}\overset{j_{C_1}}{\hookrightarrow} C_1$ be the 
generic point of $C_1$. We apply $id_{S_1}\ti {\p_2}_{S_2}$ 
on $\bar{D}$ as a cycle on $C_1\ti S_2$. 
Since the divisor $(j_{C_1}\ti id_{S_2})^*(\bar{D})$ on
$\eta_{C_1}\ti S_2$ is algebraically equivalent to a divisor 
$\eta_{C_1}\ti E$ on 
$\eta_{C_1}\ti S_2$ where $E$ is a divisor on $S_2$ defined over 
the base field $\C$, it follows that
$$(j_{C_1}\ti id_{S_2})^*(id_{S_1}\ti {\p_2}_{S_2})(\bar{D}-C_1\ti E)=(\eta_{C_1}\ti {\p_2}_{S_2})
((j_{C_1}\ti id_{S_2})^*(\bar{D})-\eta_{C_1}\ti E)=0.$$ 
So by taking the closure
of equality in $C_1\ti S_2$ it follows that
$$(id_{S_1}\ti {\p_2}_{S_2})(\bar{D})=(id_{S_1}\ti {\p_2}_{S_2})(C_1\ti E)+\sum_k p_k\ti S_2$$
for a set $\{p_k\}$ of closed points on $S_1$. In this way we 
are reduced to the case where each component of $Z$ is a 
product of two curves or is of the form $pt\ti S_2$ or $S_1\ti pt.$ 
By Conjecture (B) for surfaces one can see that the projector 
${\p_3}_{S_1}\ti {\p_2}_{S_2}$ kills the cycles of this form in 
$CH^2(X)_\Q.$

Kenichiro Kimura

Institute of Mathematics

University of Tsukuba

Tsukuba, Ibaraki

305-8571

Japan

email: kimurak@math.tsukuba.ac.jp


\begin{thebibliography}{Lampoxlong}




\bibitem[deM]{deM}del Angel, P.L., M\"uller-Stach, S.:
\newblock{Motives of uniruled $3$-folds.}
\newblock{Compositio Math.}
\newblock{\bf 112}
\newblock{(1998),}
\newblock{1-16.}


\bibitem[deM2]{deM2}del Angel, P.L., M\"uller-Stach, S.:
\newblock{On Chow motives of 3-folds.}
\newblock{Trans. Amer. Math. Soc.}
\newblock{\bf 352}
\newblock{(2000),}
\newblock{1623-1633.}

\bibitem[DM]{DM}Deninger, C., Murre, J.P.:
\newblock{Motivic decomposition of abelian schemes and the Fourier transform. }
\newblock{J. Reine Angew. Math.}
\newblock{\bf 422}
\newblock{(1991),}
\newblock{201-219.}


\bibitem[GM]{GM}Gordon, B., Murre, J.P.:
\newblock{Chow motives of elliptic modular threefolds.}
\newblock{ J. reine angew. Math.}
\newblock{\bf 514}
\newblock{(1999),}
\newblock{145--164.}

\bibitem[GHM]{GHM}Gordon, B., Hanamura, M., Murre, J.P.:
\newblock{Relative Chow-Kunneth projectors for modular varieties.}
\newblock{ J. reine angew. Math.}
\newblock{\bf 558}
\newblock{(2003),}
\newblock{1--14.}

\bibitem[GHM2]{GHM2}Gordon, B., Hanamura, M., Murre, J.P.:
\newblock{Absolute Chow-Kunneth projectors for modular varieties.}
\newblock{ J. reine angew. Math.}
\newblock{\bf 580}
\newblock{(2005),}
\newblock{139--155.}

\bibitem[Ja]{Ja}Jannsen, U.:
\newblock{Motivic sheaves and filtrations on Chow groups. }
\newblock{Proc. Sympos. Pure Math.}
\newblock{\bf 55 Part 1}
\newblock{(1994),}
\newblock{245--302.}


\bibitem[Ki]{Ki}Kimura, S.:
\newblock{Chow groups are finite dimensional, in some sense.}
\newblock{Math. Ann.}
\newblock{\bf 331}
\newblock{(2005),}
\newblock{173--201.}



\bibitem[KMP]{KMP}Kahn, B., Murre, J.P., Pedrini, C.:
\newblock{On the transcendental part of the motive of a surface.}
\newblock{preprint,}
\newblock{(2005).}



\bibitem[Mu]{Mu}Murre, J.P.:
\newblock{On the motive of an algebraic surface}
\newblock{ J. reine angew. Math.}
\newblock{\bf 409}
\newblock{(1990),}
\newblock{190-204.}

\bibitem[Mu2]{Mu2}Murre, J.P.:
\newblock{On a conjectural filtration on the Chow groups of an algebraic variety,(I and II), }
\newblock{Indag. Math. (N.S.) }
\newblock{\bf 4}
\newblock{(1993),}
\newblock{177-201.}

\bibitem[Pic]{Pic}Miller, A., M\"uller-Stach, S., Wortmann, S., Yang, Y.-H., Zuo, K.:
\newblock{Chow-K\"unneth decomposition for universal families over Picard
modular surfaces}
\newblock{preprint,}
\newblock{(2005).}





\bibitem[Sh]{Sh}Shermenev, A. M.
\newblock{Motif of an Abelian variety.}
\newblock{ Funckcional. Anal. i Prilo\v zen. }
\newblock{\bf 8}
\newblock{(1974),}
\newblock{55--61.}


\end{thebibliography}
\end{document}